\documentclass[a4paper]{amsart}
\usepackage[utf8]{inputenc}
\usepackage[T1]{fontenc}

\usepackage[style=alphabetic,maxbibnames=10,backend=bibtex]{biblatex}
\appto{\bibsetup}{\emergencystretch=0.75em}
\usepackage{paralist}
\usepackage{hyperref}

\makeatletter
\newcommand{\proofstep}[1]{%
  \par%
  \addvspace{\medskipamount}%
  \textit{#1\@addpunct{.}}\enspace\ignorespaces
}  
\newcommand{\afterlastproofstep}{%
  \par%
  \addvspace{\medskipamount}%
}  
\makeatother

\makeatother
	
\newenvironment{acknowledgement}%
  {\par%
   \hspace{\parindent}\textsc{Acknowledgements.}\enskip\ignorespaces}
  {\par}

\usepackage{amsmath, amssymb, amsfonts, mathtools, amsthm, mathrsfs, dsfont, fge, extpfeil}
\usepackage{tikz,tikz-cd,multirow}
\usetikzlibrary{matrix,arrows,positioning,calc}

\newcommand\Z{\mathbb{Z}}
\newcommand\N{\mathbb{N}}
\newcommand\M{\mathcal{M}}
\newcommand\cN{\mathcal{N}}
\newcommand\T{\mathbb{T}}
\newcommand\cT{\mathcal{T}}

\newcommand\del{\partial}
\newcommand\eps{\varepsilon}

\newcommand\cF{\mathcal{F}}
\newcommand\cR{\mathcal{R}}
\newcommand\cP{\mathcal{P}}

\newcommand\cQ{\mathcal{Q}}
\newcommand\X{\mathbb{X}}
\newcommand\cX{\mathcal{X}}
\newcommand\cY{\mathcal{Y}}

\newcommand\id{\mathrm{id}}

\DeclareMathOperator{\im}{im}
\DeclareMathOperator{\coker}{coker}

\DeclareMathOperator{\onto}{\twoheadrightarrow}

\DeclareMathOperator{\mods}{mod}

\DeclareMathOperator{\add}{add}
\DeclareMathOperator{\Hom}{\mathrm{Hom}}
\DeclareMathOperator{\End}{\mathrm{End}}
\DeclarePairedDelimiter{\Euler}{\langle}{\rangle}

\DeclareMathOperator{\dimv}{\underline{\dim}}
\DeclareMathOperator{\Tr}{Tr}
\DeclareMathOperator{\pdim}{p.dim}

\DeclareMathOperator{\gldim}{gl.dim}

\def\Ext#1#2{\mathrm{Ext}_{#1}^{#2}}
\def\defined#1{\textbf{#1}}
\newcommand{\iso}{\xrightarrow{\sim}}
\def\isom{\mathrel{\cong}}

\def\divides{\mathrel{\mid}}

\theoremstyle{definition}
\newtheorem{mydef}{Definition} 

\theoremstyle{plain}
\newtheorem{theorem}[mydef]{Theorem}
\newtheorem*{theorem*}{Theorem}
\newtheorem*{maintheorem}{Main Theorem}

\newtheorem{proposition}[mydef]{Proposition}
\newtheorem{lemma}[mydef]{Lemma}
\newtheorem{corollary}[mydef]{Corollary}
\newtheorem*{corollary*}{Corollary}

\theoremstyle{remark}

\newtheorem*{example*}{Example}

\title{Tame (hereditary) algebras of amenable representation type}
\author{Sebastian Eckert}
\address{Sebastian Eckert, Fakultät für Mathematik, Universität Bielefeld, Postfach 100~131, 33501~Bielefeld, Germany.}
\email{seckert@math.uni-bielefeld.de}
\thanks{The author has been supported~by the Alexander von~Humboldt Foundation in the framework of an Alexander von Humboldt Professorship endowed by the German Federal Ministry of Education and Research.}


\begin{document}

\begin{abstract}
We show that all finite dimensional, tame hereditary $k$-algebras are of amenable representation type (in the sense of Elek) for all fields $k$. The proof is adapted from our previous result for tame path algebras.
Further, it is proven that this results extends to some tilted algebras, in particular tame concealed algebras are amenable.
\end{abstract}

\maketitle

\section{Introduction}

The notions of hyperfiniteness for countable sets of modules and amenable representation type for algebras have been introduced by \cite{Elek2017InfiniteDimensionalRepresentationsAmenabilty}.
We will work with the definitions as follow.

\begin{mydef} \label{def:HyperfinitenessAmenability}
Let~$k$ be a field, $A$ be a finite dimensional $k$-algebra and let $\M$ be a set of finite dimensional $A$-modules. One says that $\M$ is \defined{hyperfinite} provided for every $\varepsilon > 0$ there exists a number $L_\varepsilon > 0$ such that for every $M \in \M$ there exist both, a submodule $N \subseteq M$ such that
\begin{equation} \label{eq:HFSubmoduleBig} \dim_k N \geq (1-\varepsilon) \dim_k M, \end{equation} and modules $N_1,N_2, \dots N_t \in \mods A$, with $\dim_k N_i \leq L_\varepsilon$, such that $N \isom \bigoplus_{i=1}^{t} N_i$.

The $k$-algebra $A$ is said to be of \defined{amenable representation type} provided the set of all finite dimensional $A$-modules (or more precisely, a set which meets every isomorphism class of finite dimensional $A$-modules) is hyperfinite.
\end{mydef}

Previously, the author has shown that tame path algebras of acyclic quivers are of amenable representation type (giving a new proof and not using a previous result on string algebras) while wild quiver algebras are not (using results of Elek), thus working towards \cite[Conjecture~1]{Elek2017InfiniteDimensionalRepresentationsAmenabilty}, presuming that finite dimensional algebras are of tame type if and only if they are of amenable representation type. 

In this article, we extend our previous positive result to tame hereditary algebras over arbitrary fields and to tame concealed algebras (for algebraically closed fields).

\begin{maintheorem}
Let $k$ be a field. If $A$ is a tame herdeitary, finite dimensional $k$-algebra, then $A$ is of amenable representation type.
Moreover, if $k$ is algebraically closed, any tame concealed $k$-algebra is of amenable representation type.
\end{maintheorem}

We further use the following facts from \cite{Eckert2019}.

\begin{proposition} \label{prop:AdditiveClosureStaysHyperfinite}
Let $\M$ be a family of $A$-modules. If $\M$ is hyperfinite, so is the family of all finite direct sums of modules in $\M$.
\end{proposition}

\begin{proposition} \label{prop:ExtendingHFfromSubmodulesOfBoundedCodimension} 
Let $A$ be a finite dimensional $k$-algebra. Let $\M,\cN \subset \mods A$ where $\cN$ is hyperfinite. If there is some $L \geq 0$ such that for all $M \in \M$, there exists a submodule $P \subset M$ of codimension less than or equal to $L$, and $P \in \cN$, then $\M$ is also hyperfinite.
\end{proposition}

\begin{proposition} \label{prop:HFPreservingFunctors} 
Let $k$ be a field and $A,B$ be two finite dimensional $k$-algebras. Let $F \colon \mods A \to \mods B$ be an additive, left-exact functor such that there exists $K_1, K_2 > 0$ with
\begin{equation}  \label{eq:EquivalenceCondHFPresFun}
K_1 \dim X \leq \dim F(X) \leq K_2 \dim X,\end{equation}
 for all $X \in \mods A$.
If $\cN \subseteq \mods A$ is a hyperfinite family, then the family $\{F(X) \colon X \in \cN\} \subseteq \mods B$ is also hyperfinite.
\end{proposition}

\section{Setup for tame hereditary algebras}
We will follow \cite[Section~1]{BaerGeigleLenzing1987PreprojectiveAlgebraTameHereditaryArtinAlgebra} and \cite[Section~1]{Ringel1979InfiniteDimensionalRepresentationsFDHereditaryAlgebras} and refer the reader to these papers for notions not introduced here and further references.

In the following, let $k$ be a (commutative) field and $A$ a finite dimensional $k$-algebra. We will assume that $A$ is basic, i.e. the regular module $_A A$ is is a direct sum of
pairwise non-isomorphic indecomposable projective modules.
We will focus on \emph{hereditary} algebras $A$, i.e. algebras such that all submodules of projective $A$-modules are projective again, or equivalently such that $\gldim A \leq 1$.

For $M \in \mods A$, the \emph{dimension vector of $M$} is given by
\[\dimv(M) = \left(\dim_{\End(P_i)} \Hom(P_i,M)\right)_{i=1,\dots,n},\]
where $(P_1, \dots, P_n)$ is a complete system of pairwise non-isomorphic indecomposable projective modules. We call $n$ the \emph{rank of $A$}.
Recall that the \emph{Grothendieck group $K_0(A)$} with elements $[M]$ for $A$-modules $M$ is the abelian group generated by the isomorphism classes of indecomposable objects of $\mods A$ subject to the relation that $[Y] = [X] + [Z]$ whenever $0 \to X \to Y \to Z$ is a short exact sequence in $\mods A$. Here, $K_0(A)$ coincides with $\Z^n$.

The formula
\[\Euler{\dimv M, \dimv N} \colon= \dim_k \Hom(M,N) - \dim_k \Ext{A}{1}(M,N),\]
defines a bilinear form on $K_0(A)$ and we call it the \emph{Euler form of $A$}. This homological form agrees with the bilinear form from \cite[Section~1]{Ringel1979InfiniteDimensionalRepresentationsFDHereditaryAlgebras}. We denote the associated quadratic form by $q$.

We will further follow \cite{CrawleyBoevey1991RegularModulesTameHereditaryAlgebras} in using
\begin{mydef}
A finite dimensional $k$-algebra $A$ is called \defined{tame hereditary} provided $A$ is hereditary and $A$ is of tame representation type, i.e. its quadratic form $q$ is semidefinite but not positive definite.
\end{mydef}

By $\tau(M) := D \Tr M$, we denote the $k$-dual of the Auslander-Bridger transpose and call it the \emph{Auslander-Reiten translation}. For hereditary $A$, $\tau$ is functorial, vanishes on the projective modules and induces an equivalence with inverse $\tau^{-1} = \Tr D$ from the category of all finite-dimensional $A$-modules without nontrivial projective summands into the category of all finite-dimensional $A$-modules without nontrivial injective summands. Restriction to the former subcategory gives an exact functor.
We further use the \emph{Coxeter transformation} $c : K_0(A) \to K_0(A)$, with $c([M]) = [\tau M]$ for a nonprojective indecomposable $A$-module $M$. 

Recall that using $\tau$, we can distinguish up to three kinds of indecomposable modules in $\mods A$: The \emph{preprojective} indecomposable modules are of the form $\tau^{-m}(P_i)$ for some non-negative integer $m$ and some indecomposable projective $P_i$. 
The \emph{postinjective} indecomposable modules are defined in a dual way. The remaining indecomposable modules and their finite direct sums are called \emph{regular}. 

For the moment, let us assume that $A$ is tame. The radical $\{x \in K_0(A) \colon q(x) = 0\}$ of the associated semidefinite quadratic form $q$ is then generated by a vector $\underline{h}_A$ whose entries are positive integers and at least one of them is one. We call $\underline{h}_A$ the \emph{minimal positive radical element of $A$}. Further recall that there exists an indecomposable regular module $S$ with $\dimv(S) = \underline{h}_A$. Any indecomposable $A$-module $M$ with $\dimv M$ a multiple of $\underline{h}_A$ will be called \emph{homogeneous}.
Related with $S$ is the \emph{defect of $A$}, the linear form $\Euler{\dimv S,-}$ on $K_0(A)$, which we normalize such that $\del(P_i) = -1$ for some $i \in \{1,\dots,n\}$. Note that it is invariant under $c$.
One of the main properties of the defect is that for indecomposable $M$, $\del M = \del(\dimv M) < 0$ ($\del M = 0$, $\del M > 0$) if and only if $M$ is preprojective (regular, postinjective, respectively).
The closure of preprojective, regular and postinjective modules under finite direct sums will be denoted by $\cP(A)$, $\cR(A)$ and $\cQ(A)$, respectively.


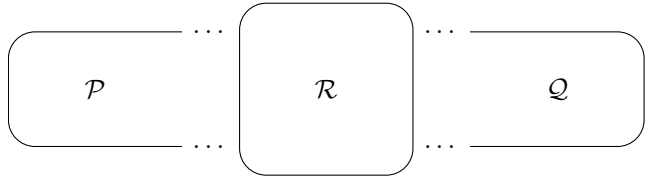
\begin{figure} \label{fig:THAlgComponents}
\begin{tikzpicture}[scale=0.38]
\begin{scope}
\clip(-11,0) rectangle (-5,4);
\draw[rounded corners=10pt] (-11,0) rectangle (-3,4);
\end{scope}
\node at (-8,2) {$\cP$};
\node at (-4,0) {\dots};
\node at (-4,4) {\dots};

\draw[rounded corners=10pt] (-3,-1) rectangle (3,5);
\node at (0,2) {$\cR$};
\node at (4,0) {\dots};
\node at (4,4) {\dots};
\begin{scope}
\clip(11,0) rectangle (5,4);
\draw[rounded corners=10pt] (3,0) rectangle (11,4);
\end{scope}
\node at (8,2) {$\cQ$};

\end{tikzpicture}
\caption{Structure of the module category $\mods A$ of a tame hereditary algebra $A$. Within in the picture, non-zero morphisms only exist from left to the right.}
\end{figure}

In the tame case, the finite-dimensional regular modules form an exact abelian subcategory $\cR = \cR(A)$. We consider its simple objects and call them \emph{regular simple}. If $M$ is in $\cR$, the sum of the regular simple submodules of $M$ is called the \emph{regular socle} of $M$, and the length $n$ of a regular composition sequence is called the \emph{regular length of $M$}. Given a regular simple module $S$ and $n \in \N$, there is a unique indecomposable regular module $S_{[n]}$ with regular socle $S$ and of regular length $n$, and every indecomposable regular module is of this form.
Recall that $\tau$ operates on the set of (isomorphism classes of) regular simple modules with finite orbits, and all but at most three orbits are one element sets. Let $\X$ be the set of these orbits. If $S$ and $S'$ are regular simple, then $\Ext{A}{1}(S,S') \neq 0$ if and only if $S' = \tau S$. Thus, the category $\cR$ decomposes as the direct sum of uniserial categories $R_t(A)$ with finitely many simple objects, where $t$ runs through the set $\X$:
\[\cR(A) \isom \coprod_{t \in \X} R_t(A).\]
We may call the poset formed by the indecomposable objects in such a uniserial category $R_t$ a \emph{tube} and use $\T$ to denote it for a fixed $t \in \X$. We say that the regular-simple modules form the \emph{mouth} of the tube $\T$.

For $S$ regular simple in $R_t$, let $n_t$ denote the \emph{rank} of the corresponding tube, i.e. the smallest positive integer with $\tau^{n_t} \isom S$. Note that $S^{n_t}$ is always homogeneous, whereas modules $S^i$, with $1 \leq i < n_t$ are not homogeneous. We call $R_t$ homogeneous iff all regular simple modules in $R_t$ are.

Given a $k$-algebra $A$ and an $A$-bimodule $M$, on which $k$ acts centrally, we can form the \emph{tensor algebra}
\[T_A(M) := \bigoplus_{n\geq 0} M^{\otimes n} = A \oplus M \oplus (M\otimes_A M) \oplus (M \otimes_A M \otimes_A) \oplus \dots .\]
If $A$ is semisimple, $T_A(M)$ will be hereditary.

Recall the following notions from \cite[Section~1]{DlabRingel1978RepresentationsTameHereditaryAlgebras}: Given a field $k$, a $k$-automorphism $\eps$ and an $(\eps,1)$-derivation $\delta$ of $k$, we can define the $k$-$k$-bimodule $M(\eps,\delta)$ which as a left $k$ vector space is $_k k \oplus _k k$, while the right $k$-action involves the derivation $\delta$. 
Using this information, we denote by $\tilde{A}_n(\eps,\delta)$ the $(n+1)\times(n+1)$-matrix ring
\[\begin{pmatrix} k & k & \dots & k & M(\eps,\delta) \\ & k & \dots & k & k \\ & & \ddots & \vdots & \vdots\\ & 0 & & k & k \\ & & & & k \end{pmatrix}.\]

\begin{theorem} \cite[Corollary~2]{DlabRingel1978RepresentationsTameHereditaryAlgebras}
A hereditary finite  dimensional $k$-algebra $A$ is of tame representation type if and only if it is Morita equivalent to the product of a tame tensor algebra and a finite number of algebras of the form $\tilde{A}_n(\eps,\delta)$.
\end{theorem}


\section{The proof}
In this section, we shall develop the proof of the Main Theorem.
A first lemma pertaining the existence of preprojective submodules of bounded codimension of regular modules will be used in the next subsection but also later on.

\begin{lemma} \label{lemma:PPSubmodulesOfRegulars}
Let $A$ be a finite dimensional tame hereditary algebra.
Denote by $l$ the maximum of the dimensions of the indecomposable injective modules.
Let $R$ be an indecomposable regular module. Then there exists a preprojective submodule $U\subseteq R$ such that \[\dim_k R - \dim_k U \leq l.\]
\begin{proof}
Let $S$ be the regular socle of $R$. Then there exists a non-zero map $f \colon S \to I$, where $I$ is an indecomposable injective module. 
Since $I$ is injective and $\iota \colon S \to R$ is injective, this map lifts to a map $\bar{f} \colon R \to I$. 
Denote $U := \ker \bar{f}$ and $K := \ker f$. Using the Snake Lemma, we get a exact commutative diagram
\begin{center}
\begin{tikzpicture}
\matrix (m) [matrix of math nodes,nodes={anchor=center},row sep=2em,column sep=3em]
  {
	 & 0 & 0 \\
     0 & K & U & R/S & 0\\
     0 & S & R & R/S & 0\\
     0 & I & I & 0\\
     & 0 & 0\\
     };
  \path[->]
    (m-1-2) edge (m-2-2)
    (m-1-3) edge (m-2-3)
    (m-2-1) edge (m-2-2)
    (m-2-2) edge (m-2-3)
			edge (m-3-2)
	(m-2-3) edge (m-2-4)
			edge (m-3-3)
	(m-2-4) edge (m-2-5)
			edge [-,double equal sign distance] (m-3-4)
	(m-3-1) edge (m-3-2)
	(m-3-2) edge [right hook->] node [above] {$\iota$} (m-3-3)
			edge node [left] {$f$} (m-4-2)
	(m-3-3) edge (m-3-4)
			edge node [left] {$\bar{f}$} (m-4-3)
	(m-3-4) edge (m-3-5)
			edge (m-4-4)
	(m-4-1) edge (m-4-2)
	(m-4-2) edge [-,double equal sign distance] (m-4-3)
			edge (m-5-2)
	(m-4-3) edge (m-4-4)
			edge (m-5-3);
\end{tikzpicture}
\end{center}
Clearly, $K$ and $U$ as subobjects of regular modules cannot have postinjective summands. 
Assume that $U$ had a regular summand. Then this summand must contain $S$ resp. $\iota(S)$, for this is the smallest regular submodule of $R$, using the fact that the regular modules form a uniserial category and no other tubes map to $R$. But then we have
\[ 0 = \bar{f}(U) \supset \bar{f}(\iota S) = f(S) \neq 0,\]
a contradiction. Hence $U$ can only have preprojective summands and is therefore preprojective.

\end{proof}
\end{lemma}

\subsection{Finding large submodules for rank two algebras}
We will first restrict to the case of tame hereditary algebras of rank two, the non-algebraically closed analogue of the 2-Kronecker quiver case. As $\tilde{A}_1(\eps,\delta)$ is a tensor algebra, the representation theory of all tame hereditary algebras of rank two is covered by the methods of \cite{DlabRingel1976IndecomposableRepresentationsGraphsAlgebras}.

\begin{lemma} \label{lemma:DescentOnDefectOfPIs} 
Let $X$ be some indecomposable postinjective $kQ$-module of defect $\del(X) = d$. Then there is an injective module $I(j)$ such that there exists a non-zero morphism $\theta \colon X \to I(j)$. Moreover,
for any direct summand $Z$ of $\ker \theta$, we have $\del(Z) < d$.
\begin{proof}
Let $E(X)$ be the injective envelope of $X$, and take $I(j)$ to be some direct summand of $E(X)$. This yields a non-zero homomorphism $\theta \colon X \to E(X) \to I(j)$.
Consider the exact sequence
$$0 \to \ker \theta \to X \to \im \theta \to 0.$$
Since there is a map from a postinjective module to $\im \theta$, the latter must be postinjective or zero. Yet, $\im \theta \neq 0$, since $\theta$ is non-zero. Thus, $\del(\im \theta) > 0$. 
Therefore, $\del(\ker \theta) < \del(X)$.
If $\ker \theta$ only had preprojective or regular summands $Z$, we are done, for then $\del(Z) \leq 0$.
Thus, we may assume that there is some some postinjective direct summand $Z$. Since $Z$ embeds into $\ker \theta$ and the kernel embeds into $X$, we get a short exact sequence
\[0 \to Z \to X \to X/Z \to 0.\]
Since $X$ is postinjective, again $X/Z$ must be postinjective or zero. If it was zero, then $Z \isom \ker \theta \isom X$, a contradiction, since $\del(\ker \theta) \neq \del(X)$. Thus, $\del(X/Z) > 0$, and hence we may conclude
$\del(Z) < \del(X) = d$.
\end{proof}
\end{lemma}

\begin{lemma} \label{lemma:CoxeterRank2ChangeByhQ}
Let $A$ be a tame hereditary algebra of rank two. 
Then there is $g \in \N^{+}$ such that 
\[c^{\pm}(x) = x \pm g \underline{h}_A \del(x), \quad \text{ for all } x \in K_0(A).\]
\begin{proof}
Recall that  the defect is given by $\del(x) = \frac{1}{r} \Euler{\underline{h}_A,x}$ for some $r \in \N^{+}$.
We show the claim by a case by case analysis. There are two types of algebras of rank two, namely type $(2,2)$ and $(1,4)$. For the latter, we also take into consideration the two possible orientations.
\proofstep{type $(2,2)$}
Here, we have $\underline{h}_A = (1,1)$ and the Coxeter transformation (resp. its inverse) as a matrix is given by
\[c = \begin{pmatrix} -1 & 2\\-2 & 3\end{pmatrix}.\]
Now \begin{align*} c(x_1,x_2) = (-x_1+2x_2,-2x_1+3x_2) &= (x_1,x_2) + (-2x_1+2x_2,-2x_1+2x_2)\\ &= (x_1,x_2) + (2,2) (-x_1+x_2).\end{align*}
On the other hand, we have
\begin{align*} \del(x) & = \Euler{\underline{h}_A,x} = \Euler{e_1,x} + \Euler{e_2,x}\\
&= \left(x_1 \dim_k A_1 - x_2 \dim_k \eps_1 M \eps_2\right) + \left( x_2 \dim_k A_2 - x_1 \dim_k \eps_2 M \eps_1 \right)\\
&= x_1 + x_2 - 2x_1 = -x_1+x_2. \end{align*}
\proofstep{type $(1,4)$}
Here, we have $\underline{h}_A = (2,1)$ and the Coxeter transformation (resp. its inverse) as a matrix is given by
\[c = \begin{pmatrix} -1 & 4\\-1 & 3\end{pmatrix}.\]
Now \begin{align*} c(x_1,x_2) = (-x_1 + 4x_2,-x_1 + 3x_2) &= (x_1,x_2) + (-2x_1+4x_2,-x_1+2x_2)\\ &=(x_1,x_2) + (2,1) (-x_1+2x_2).\end{align*}
On the other hand, we have
\begin{align*} \del(x) & = \frac{1}{2} \Euler{\underline{h}_A,x} = \frac{1}{2} \left( 2\Euler{e_1,x} + \Euler{e_2,x}\right)\\
&= \left(x_1 \right) + \frac{1}{2} \left( 4 x_2 - 4 x_1\right)
= -x_1 + 2x_2. \end{align*}

\proofstep{type $(4,1)$}
Here, we have $\underline{h}_A = (1,2)$ and the Coxeter transformation (resp. its inverse) as a matrix is given by
\[c = \begin{pmatrix} -1 & 1\\-4 & 3\end{pmatrix}.\]
Now \begin{align*} c(x_1,x_2) = (-x_1 + x_2,-4x_1 + 3x_2) &= (x_1,x_2) + (-2x_1+x_2,-4x_1+2x_2)\\ &=(x_1,x_2) + (1,2) (-2x_1+x_2).\end{align*}
On the other hand, similar to the other orientation, we have
$\del(x) = -2x_1+x_2$.
\end{proof}
\end{lemma}

\begin{proposition} \label{prop:EpiToRegularForPPDefect-1}
Let $P$ be an indecomposable preprojective module of defect $\del(P) = -1$.
For $i\geq 0$, in type $\tilde{A}_{11}$, consider the modules $P[i] := \tau^{-i} P$ and in type $\tilde{A}_{12}$ consider the indecomposable preprojective modules with $\dimv P[i] = \dimv P + i \underline{h}_A$. Choose some non-zero homomorphism $P \to P[1]$.
Let $R$ be the regular module given as the cokernel of this map and denote by $R_{[m]}$ the unique indecomposable regular module with regular-socle $R$ and of regular-length $m$.
Then for any $n > 0$ and $m \leq n$, there exists an epimorphism
$\phi_{n,m} \colon P[n] \onto R_{[m]}$ with kernel $P[n-m]$.
\begin{proof}
Let $P$ be an indecomposable projective module of defect $\del(P) = -1$. For type $\tilde{A}_{11}$, let $P[1] := \tau^{-1}P$ be its inverse Auslander Reiten-translate. For type $\tilde{A}_{12}$, let $P[1]$ be the uniquely determined indecomposable preprojective module with dimension vector $\dim P + \underline{h}_A$.
Then $\Euler{P,P[1]} >0$, so there is a non-zero map $\iota \colon P \to P[1]$. Since $\del(P) = -1$, this map must be injective. Consider the short exact sequence
\[ \eta \colon 0 \to P \xrightarrow{\iota} P[1] \xrightarrow{\pi} R \to 0, \]
where $R = \coker \iota$. 
Since $\del(P) = \del(P[1]) = -1$, we have $\del(R) = 0$.
If $R$ had a postinjective summand $Q$ of defect $\geq 1$, we get the pullback
\begin{center}
\begin{tikzpicture}
\matrix (m) [matrix of math nodes,nodes={anchor=center},row sep=2em,column sep=3em]
  {
     0 & P & T & Q & 0\\
     0 & P & P[1] & R & 0,\\};
  \path[->]
    (m-1-1) edge node [above] {} (m-1-2)
    (m-1-2) edge [right hook->] node [above] {} (m-1-3)
			edge [-,double equal sign distance] node [above] {} (m-2-2)
    (m-1-3) edge [->>] node [above] {} (m-1-4)
			edge [right hook->] node [right] {} (m-2-3)
    (m-1-4) edge node [above] {} (m-1-5)
			edge [right hook->] node [right] {} (m-2-4)
    (m-2-1) edge node [above] {} (m-2-2)
    (m-2-2) edge [right hook->] node [above] {$\iota$} (m-2-3)
    (m-2-3) edge [->>] node [above] {$\pi$} (m-2-4)
    (m-2-4) edge node [above] {} (m-2-5);
\end{tikzpicture}
\end{center}
where $\del(T) = \del(P) + \del(Q) \geq 0$, but $T$ is also a submodule of $P[1]$, hence preprojective. This implies that $T=0$.
Thus, $R$ cannot have postinjective and hence neither preprojective summands. It is therefore regular.
By construction, as $\del(P) = -1$ and due to Lemma~\ref{lemma:CoxeterRank2ChangeByhQ}, we have that $\dimv R = h_A$. As every (indecomposable) submodule of $R$ would have dimension vector belonging to an indecomposable preprojective or preinjective module, we see that $R$ must be simple regular.

Next, consider the regular module $R_{[2]}$ such that
\[\xi \colon 0 \to R \to R_{[2]} \to R \to 0,\]
is an almost split sequence, i.e. $0 \neq \xi \in \Ext{A}{1}(R,R)$. By applying $\Hom_A(R,-)$ to $\eta$, we get a long exact sequence
\[0 \to (R,P) \xrightarrow{\iota_*} (R,P[1]) \xrightarrow{\pi_*} (R,R) \xrightarrow{\eta_*} {^1(R,P)} \xrightarrow{\iota_*} {^1 (R,P[1])} \xrightarrow{\pi_*} {^1 (R,R)} \to 0,\]
as $A$ is hereditary.
This shows that the map $\pi_* \colon \Ext{A}{1}(R,P[1]) \to \Ext{A}{1}(R,R)$ is surjective. Now, there exists $\zeta \in \Ext{A}{1}(R,P[1])$ s.t. $\pi_*(\zeta) = \xi$, where $\pi_*$ is the push-out map.
Hence
\begin{center}
\begin{tikzpicture}
\matrix (m) [matrix of math nodes,nodes={anchor=center},row sep=2em,column sep=3em]
  {
	 & 0 & P & \ker \phi & 0 \\
     \zeta\colon & 0 & P[1] & E & R & 0\\
     \xi = \pi_*(\zeta) \colon & 0 & R & R_{[2]} & R & 0\\
     & & \coker \pi & \coker \pi' & 0\\
     };
  \path[->]
    (m-1-2) edge (m-1-3)
    (m-1-3) edge (m-1-4)
			edge [right hook->] node [left] {$\iota$} (m-2-3)
	(m-1-4) edge (m-1-5)
			edge [right hook->] (m-2-4)
	(m-1-5) edge (m-2-5)
    (m-2-2) edge (m-2-3)
	(m-2-3) edge [right hook->] node [above] {$\sigma$} (m-2-4)
			edge [->>] node [left] {$\pi$} (m-3-3)
	(m-2-4) edge [->>] node [above] {$\theta$} (m-2-5)
			edge [->>] node [left] {$\phi$} (m-3-4)
	(m-2-5) edge (m-2-6)
			edge [-,double equal sign distance] (m-3-5)
	(m-3-2) edge (m-3-3)
	(m-3-3) edge [right hook->] node [above] {$\alpha$} (m-3-4)
			edge (m-4-3)
	(m-3-4) edge [->>] node [above] {$\beta$} (m-3-5)
			edge (m-4-4)
	(m-3-5) edge (m-3-6)
			edge (m-4-5)
	(m-4-3) edge (m-4-4)
	(m-4-4) edge (m-4-5);
  \draw[->,red,rounded corners] 
	(m-1-5)	-| ($(m-2-6.east)+(.5,0)$) |- ($(m-3-4)!0.35!(m-2-4)$) -| ($(m-3-2.west)+(-.5,0)$) |- (m-4-3);
\end{tikzpicture}
\end{center}
is an exact commutative diagram, using the Snake Lemma.
Alltogether, we get the exact commutative diagram
\begin{center}
\begin{tikzpicture}
\matrix (m) [matrix of math nodes,nodes={anchor=center},row sep=2em,column sep=3em]
  {
	 & & & 0 & 0\\
	 \eta \colon & 0 & P & P[1] & R & 0 \\
     & 0 & \ker \phi & E & R_{[2]} & 0\\
     & & & R & R\\
     & & & 0 & 0\\
     };
  \path[->]
    (m-1-4) edge (m-2-4)
    (m-1-5) edge (m-2-5)
    (m-2-2) edge (m-2-3)
	(m-2-3) edge [right hook->] node [above] {$\iota$} (m-2-4)
			edge node [left] {$\isom$} (m-3-3)
	(m-2-4) edge [->>] node [above] {$\pi$} (m-2-5)
			edge [right hook->] node [left] {$\sigma$} (m-3-4)
	(m-2-5) edge (m-2-6)
			edge [right hook->] node [left] {$\alpha$} (m-3-5)
	(m-3-2) edge (m-3-3)
	(m-3-3) edge [right hook->] node [above] {$\gamma$} (m-3-4)
	(m-3-4) edge [->>] node [above] {$\phi$} (m-3-5)
			edge [->>] node [left] {$\theta$} (m-4-4)
	(m-3-5) edge (m-3-6)
			edge [->>] node [left] {$\beta$} (m-4-5)
	(m-4-4) edge [-,double equal sign distance] (m-4-5)
			edge (m-5-4)
	(m-4-5) edge (m-5-5);
\end{tikzpicture}
\end{center}
where $\ker \phi \isom P$ via the Snake Lemma, for $\ker(R \iso R) = 0$.

Now assume that $E$ has a regular summand, e.g. $E = R'\oplus E'$.
Since the regular component is uniserial, we must have
\begin{enumerate}[i)]
\item $\phi(R') \subseteq \alpha(R_{})$ or
\item $\alpha(R_{}) \subset \phi(R') \subseteq R_{[2]}$.
\end{enumerate}
Since $\alpha(R_{})$ is a maximal submodule of $R_{[2]}$, ii) implies $\phi(R') = R_{[2]}$.
\proofstep{i)} Applying $\beta$ yields $\theta(R') = \beta \phi(R') \subseteq \beta \alpha(R_{})=0$, so $R' \subseteq \ker \theta \isom \im \sigma \isom P[1]$, thus $R'$ is preprojective for $A$ is hereditary. A contradiction.
\proofstep{ii)} This implies that $\phi_{|R'} \colon R' \to R_{[2]}$ is surjective.
On the other hand, $\ker(\phi_{|R'}) = \ker (\phi) \cap R' \subseteq \ker(\phi) \isom P$, so is preprojective. But $\phi_{|R'}$ is a map between regular modules, thus has a regular kernel. Hence $\ker(\phi_{|R'}) = 0$, and $\phi_{|R'} \colon R' \to R_{[2]}$ is bijective. But then there exists a section $\phi' \colon R_{[2]} \to R' \subset E$ such that
$\phi \circ \phi' = \id_{R_{[2]}}$. This is equivalent to the existence of a retraction $\gamma' \colon E \to P$ such that $\gamma' \circ \gamma = \id_P$. But then $\gamma' \circ \sigma \circ \iota = \gamma' \circ \gamma = \id_P$, so $\eta$ splits. A contradiction.

Alltogether, $E$ must be preprojective. Note that as an extension of preprojective and regular modules it cannot be postinjective.
Moreover, $\del(E) = \del(P[1]) + \del(R) = \del(P[1]) = -1$, so it is indecomposable. Besides,
\[\dimv E = \dimv P[1] + \dimv R = \dimv P[1] + (\dimv P[1] - \dimv P) = \dimv P[2].\]
Hence $E \isom P[2]$, since both are rigid modules and thus determined by their dimension vectors.
We have therefore constructed a surjective map $\phi \colon P[2] \to R_{[2]}$ with kernel $P$.

Now assume we have constructed a surjective map $\pi_m \colon P[m] \to R_{[m]}$ for some $m \in \N$ with kernel $P$. Consider the exact sequence
\[\eta_m \colon 0 \to P \xrightarrow{\iota_m} P[m] \xrightarrow{\pi_m} R_{[m]} \to 0.\]
As in the base case, the pushout map $\pi_* \colon \Ext{A}{1}(R,P[m]) \to \Ext{A}{1}(R,R_{[m]})$ is surjective. Thus, given the standard AR-sequence
\[\xi_m \colon 0 \to R_{[m]} \xrightarrow{\alpha_m} R_{[m+1]} \xrightarrow{\beta_m} R \to 0,\]
there exists some $\zeta_m \in \Ext{A}{1}(R,P[m])$,
\[\zeta_m \colon 0 \to P[m] \xrightarrow{\sigma_m} E_m \xrightarrow{\theta_m} R \to 0,\]
such that we get the exact commutative diagram
\begin{center}
\begin{tikzpicture}
\matrix (m) [matrix of math nodes,text height=1.5ex, text depth=0.25ex,,row sep=2em,column sep=3em]
  {
	 & & & 0 & 0\\
	 \eta_m \colon & 0 & P & P[m] & R_{[m]} & 0 \\
     & 0 & \ker \phi_m & E_m & R_{[m+1]} & 0\\
     & & & R & R\\
     & & & 0 & 0\\
     };
  \path[->]
    (m-1-4) edge (m-2-4)
    (m-1-5) edge (m-2-5)
    (m-2-2) edge (m-2-3)
	(m-2-3) edge [right hook->] node [above] {$\iota_m$} (m-2-4)
			edge node [left] {$\isom$} (m-3-3)
	(m-2-4) edge [->>] node [above] {$\pi_m$} (m-2-5)
			edge [right hook->] node [left] {$\sigma_m$} (m-3-4)
	(m-2-5) edge (m-2-6)
			edge [right hook->] node [left] {$\alpha_m$} (m-3-5)
	(m-3-2) edge (m-3-3)
	(m-3-3) edge [right hook->] node [above] {$\gamma_m$} (m-3-4)
	(m-3-4) edge [->>] node [above] {$\phi_m$} (m-3-5)
			edge [->>] node [left] {$\theta_m$} (m-4-4)
	(m-3-5) edge (m-3-6)
			edge [->>] node [left] {$\beta_m$} (m-4-5)
	(m-4-4) edge [-,double equal sign distance] (m-4-5)
			edge (m-5-4)
	(m-4-5) edge (m-5-5);
\end{tikzpicture}
\end{center}
As in the base case, $\ker \phi_m \isom P$.
Assume that $E_m$ has a regular summand $R'$, i.e. $E_m \isom R' \oplus E'$.
Since the regular component is uniserial, we must have
\begin{enumerate}[i)]
\item $\phi_m(R') \subseteq \alpha_m(R_{[m]})$ or
\item $\alpha_m(R_{[m]}) \subset \phi_m(R') \subseteq R_{[m+1]}$.
\end{enumerate}
Since $\alpha_m(R_{[m]})$ is a maximal submodule of $R_{[m+1]}$, ii) implies $\phi_m(R') = R_{[m+1]}$.
\proofstep{i)} Applying $\beta_m$ yields $\theta_m(R') = \beta_m \phi_m(R') \subseteq \beta_m \alpha_m(R_{[m]})=0$, so $R' \subseteq \ker \theta_m \isom \im \sigma_m \isom P[m]$, thus $R'$ is preprojective for $A$ is hereditary. A contradiction.
\proofstep{ii)} This implies that $\phi_m{_{|R'}} \colon R' \to R_{[m+1]}$ is surjective.
On the other hand,
\[\ker(\phi_m{_{|R'}}) = \ker (\phi_m) \cap R' \subseteq \ker(\phi_m) \isom P,\] so is preprojective. But $\phi_m{_{|R'}}$ is a map between regular modules, thus has a regular kernel. Hence $\ker(\phi_m{_{|R'}}) = 0$, and $\phi_m{_{|R'}} \colon R' \to R_{[m+1]}$ is bijective. But then there exists a section $\phi_m' \colon R_{[m+1]} \to R' \subset E_m$ such that
$\phi_m \circ \phi_m' = \id_{R_{[m+1]}}$. This is equivalent to the existence of a retraction $\gamma_m' \colon E_m \to P$ such that $\gamma_m' \circ \gamma_m = \id_P$. But then $\gamma_m' \circ \sigma_m \circ \iota_m = \gamma_m' \circ \gamma_m = \id_P$, so $\eta_m$ splits. A contradiction.

All in all, $E_m$ must be preprojective. Note that it cannot be postinjective as an extension of preprojective and regular modules.
Moreover, $\del(E_m) = \del(P[m]) + \del(R) = \del(P[m]) = -1$, so it is indecomposable. Besides,
\[\dimv E_m = \dimv P[m] + \dimv R = \dimv P[m] + (\dimv P[1] - \dimv P) = \dimv P[m+1].\]
Hence $E \isom P[m+1]$.
We have therefore constructed a surjective map $\phi \colon P[m+1] \to R_{[m+1]}$ with kernel $P$.

Thus, we may assume that there exists a surjective map $\pi_n \colon P[n] \to R_{[n]}$ for all $n>0$, which has kernel $P$. Since there is also a surjective map $R_{[n]} \to R_{[m]}$ for all $m\leq n$ by the uniseriality, we have established the existence of a surjective homomorphism $\pi_{n,m} \colon P[n] \to R_{[m]}$.
This morphism fits into a short exact sequence
\[0 \to \ker \pi_{n,m} \hookrightarrow P[n] \twoheadrightarrow R_{[m]} \to 0.\]
Since $P[n]$ is preprojective, so is $\ker \pi_{n,m}$. Since $\del(\ker \pi_{n,m}) = \del(P_n) - \del(R_{[m]}) = -1$, it is indecomposable.
Moreover,
\begin{align*}
\dimv \ker \pi_{n,m} &= \dimv P[n] - \dimv R_{[m]}\\
& = \dimv P[m] + (m-n) \left(\dimv P[1] - \dimv P\right) - (\dimv P[m] - \dimv P)\\
& = \dimv P[n-m],\end{align*}
using the recursion formula $\dimv P[i+1] = \dimv P[i] + (\dimv P[1] - \dimv P)$ along with $\dimv R_{[m]} = \dimv P[m] - \dimv P$.
Hence $\ker \pi_{n,m} \isom P[n-m]$.
\end{proof}
\end{proposition}

\subsubsection*{Algorithm for preprojectives of defect $-1$ for rank two}

Let $A$ be a tame hereditary algebra of rank two, with minimal radical vector $\underline{h}_A$ and $h_A := \sum_{j=1}^{n} f_j {\underline{h}_A}_i$ the $k$-weighted sum of the entries.
Let $\underline{q}$ be the dimension vector of the indecomposable injective $Q$ of maximal $k$-dimension, let $q = \sum_{j=1}^{n} \dim_k Q_j$.
Let $1 > \varepsilon > 0$. Choose $L_{\varepsilon} := \max\{\frac{2q}{\varepsilon},2 g h_A\}$ where $g$ is from Lemma~\ref{lemma:CoxeterRank2ChangeByhQ}.

Now let $X=P[{i_0}]$ be any indecomposable preprojective module of defect $\del(X) = -1$, where $P$ is an indecomposable projective of defect $-1$. Denote the corresponding regular modules as in Proposition~\ref{prop:EpiToRegularForPPDefect-1} by $R_{[i]}$ respectively.
If $\dim_k X \leq L_{\varepsilon}$, choosing $Y=X$ we have found a suitable submodule to prove hyperfiniteness. Hence assume that $\dim_k X > L_{\varepsilon}$.

We will now give an iterative construction involving $i_j$, $T_j = P[i_j]$ and $t_{j+1}$. We start with $j=0$.

For each $j$ we proceed as follows:
Let $R_{[t_{j+1}]}$ be a regular indecomposable such that
\[L_{\varepsilon} - gh_A \leq \dim_k R_{[t_{j+1}]} \leq L_{\varepsilon}, \quad t_{j+1} \leq {i_{j}}.\]
This is possible, since by construction $\dimv R_{[1]} = gh_A.$
Then by Propositon~\ref{prop:EpiToRegularForPPDefect-1}, we get a short exact sequence
\[0 \to P[{i_{j}}-t_{j+1}] \to P[{i_{j}}] \to R_{[t_{j+1}]} \to 0.\]
By Lemma~\ref{lemma:PPSubmodulesOfRegulars}, there exists a preprojective submodule $U_{j+1} \subset R_{[t_{j+1}]}$ such that \[\dimv R_{[t_{j+1}]} - \dimv U_{j+1} \leq \underline{q}.\]
Consider the pullback
\begin{center}
\begin{tikzpicture}
\matrix (m) [matrix of math nodes, text height=1.5ex, text depth=0.25ex, row sep=2em,column sep=3em]
  {
     0 & P[{i_{j}}-t_{j+1}] & E_{j+1} & U_{j+1} & 0\\
     0 & P[{i_{j}}-t_{j+1}] & T_{j} & R_{[t_{j+1}]} & 0.\\};
  \path[->]
    (m-1-1) edge node [above] {} (m-1-2)
    (m-1-2) edge [right hook->] node [above] {} (m-1-3)
			edge [-,double equal sign distance] node [above] {} (m-2-2)
    (m-1-3) edge [->>] node [above] {} (m-1-4)
			edge [right hook->] node [right] {} (m-2-3)
    (m-1-4) edge node [above] {} (m-1-5)
			edge [right hook->] node [right] {} (m-2-4)
    (m-2-1) edge node [above] {} (m-2-2)
    (m-2-2) edge [right hook->] node [above] {} (m-2-3)
    (m-2-3) edge [->>] node [above] {} (m-2-4)
    (m-2-4) edge node [above] {} (m-2-5);
\end{tikzpicture}
\end{center}
If the upper sequence is non-split, we have \[0 \neq \Ext{A}{1}(S,P[i_{j}-t_{j+1}]) \isom D\Hom_A(P[i_{j}-t_{j+1}],\tau S),\] for some indecomposable direct summand $S \divides U_{j+1}$ and a non-zero homomorphism $P[i_{j}-t_{j+1}] \to \tau S$, which must be a monomorphism, for $\del(P[i_{j}-t_{j+1}]) = -1.$ Hence, using Lemma~\ref{lemma:CoxeterRank2ChangeByhQ},
\begin{equation} \label{eq:AlgPPDef-1Rk2KernelLeps}
\dim P[i_{j}-t_{j+1}] \leq \dim \tau S \leq \dim S \leq L_{\varepsilon}.\end{equation}
But now the iteration terminates with $N=j+1$. By choosing $Y := E_N \oplus \bigoplus_{\ell=1}^{N-1} U_{\ell}$, we have found a suitable submodule with its summands' dimensions bound by $2 L_\varepsilon$ due to \eqref{eq:AlgPPDef-1Rk2KernelLeps} and the fact that all $R_{[t_\ell]}$ have bounded dimension, still big enough using $\dimv X = \dimv E_N + \sum_{\ell=1}^{N-1} \dimv U_{\ell} + \sum_{\ell=1}^N e_\ell$, where $e_\ell$ is an error vector bounded by $\underline{q}$.

Thus we may assume that the upper sequence splits and $E_{j+1} \isom P[i_{j}-t_{j+1}] \oplus U_{j+1}$.
Now let $i_{j+1} := i_{j}-t_{j+1}$ and proceed with step $j+1$. 
This process terminates since the dimension of $T_j$ compared to $T_{j-1}$ decreases in each ``splitting'' step by at least $L_\varepsilon - gh_A$ until it is smaller than $L_{\varepsilon}$. Note that the number of steps $N$ is bounded by
\[N \leq \frac{\dim_k X }{L_\varepsilon-g h_A} \leq \frac{\dim_k X}{\left(\frac{q}{\varepsilon} + gh_A\right) - g h_A} = \frac{\dim_k X \varepsilon}{q}.\]

Moreover, in each step we have $\dimv T_j = \dimv T_{j+1} + \dimv R_{[t_{j+1}]}$ and $\dimv R_{[t_j]} = \dimv U_j + e_j$, where $e_j$ is an error vector bounded by $\underline{q}$. Combining this in a telescope sum yields
\begin{equation} \label{eq:TelescopeSumPPDefect-1}
\dimv X = \dimv T_0 = \dots = \dimv T_N + \sum_{\ell=1}^N \left( \dimv U_\ell + e_\ell \right).
\end{equation}
Define the submodule $Y \subset X$ to be
\[Y := T_N \oplus \bigoplus_{\ell=1}^N U_\ell.\]
Here, each summand has $k$-dimension bounded by $L_{\varepsilon}$ by construction of $N$ resp. $R_{[t_\ell]}$.
By \eqref{eq:TelescopeSumPPDefect-1}, we have
\begin{align*} \dimv Y &= \dimv T_0 - \sum_{\ell=1}^N e_\ell \geq \dimv T_0 - N \cdot \underline{q} \geq \dimv M - \frac{\dim_k X \varepsilon}{q} \underline{q},
\end{align*}
implying that $\dim_k Y \geq (1-\varepsilon) \dim_k X$.

This construction thus shows the following
\begin{lemma} \label{lemma:PPDefect-1HF}
Let $A$ be a tame hereditary algebra of rank two. Let $\eps > 0$.
Then there exists some $L_\eps > 0$ such that for all indecomposable preprojective modules $X$ of defect $-1$ there exists a submodule $Y$ such that $\dim_k Y \geq (1-\eps) \dim_k X$ and each indecomposable direct summand of $Y$ has dimension bounded by $L_\eps$.
\end{lemma}

To now deal with the remaining indecomposable preprojective modules, we use the following well known Theorem due to Auslander and Reiten.
\begin{theorem}\cite[Theorem~V.5.3]{AuslanderReitenSmaloe1995RepresentationTheoryArtinAlgebras},\cite[Theorem~IV.1.10]{AssemSimsonSkowronski2010TechniquesRT} \label{thm:IrreducibleFromSinkMap}
Let $C$ be an indecomposable module and let $h \colon B \to C$ be a minimal right almost split morphism (viz. sink map). If $E \neq 0$ is some direct summand of $B$, the induced map $g \colon E \to C$ with $g = h_{|E}$ is irreducible.
\begin{proof}
We fix some notation denoting the decomposition $B = E \oplus E'$ and defining $g' = h_{|E'}$, i.e. $h = \begin{psmallmatrix}g & g'\end{psmallmatrix}$.

First we assume that $g$ is a split monomorphism, i.e. there exists some $f \colon C \to E$ such that $f \circ g = \id_{E}$ and $C \isom E \oplus \ker f$. Now, $C$ is indecomposable and $E \neq 0$, so $\ker f = 0$ and $f$ must be a monomorphism. But $f$ is also surjective. Thus $f$ is an isomorphism, and so must be $g$. By considering  $\begin{psmallmatrix}g & g'\end{psmallmatrix} \circ \begin{psmallmatrix}g^{-1}\\0\end{psmallmatrix} = \id_C$ we see that $h$ is a split epimorphism, a contradiction.

Second, assume that $g$ is a split epimorphism, i.e. there exists some $f' \colon C \to E$ such that $g \circ f' = \id_C$. Then $\begin{psmallmatrix}g & g'\end{psmallmatrix} \circ \begin{psmallmatrix}f'\\0\end{psmallmatrix} = \id_C$, which also implies that $h$ is a split epimorphism.

Now assume that $g = s \circ t$, where $t \colon E \to X$ and $s \colon X \to C$ and assume that $s$ is not a split epimorphism, again a contradiction.
Then, since $h$ is right almost split, there exists some $\eta = \begin{psmallmatrix}u \\ v \end{psmallmatrix} \colon X \to B = E \oplus E'$ with $s = h \circ \eta = \begin{psmallmatrix}g & g'\end{psmallmatrix} \circ \begin{psmallmatrix}u\\v\end{psmallmatrix}$.
We obtain the commutative diagram
\begin{center}
\begin{tikzpicture}
\matrix (m) [matrix of math nodes, text height=1.5ex, text depth=0.25ex, row sep=3em,column sep=3em]
  {
     E \oplus E' & X \oplus E' & E \oplus E'\\
     & C\\};
  \path[->]
    (m-1-1) edge node [above] {$\begin{psmallmatrix}t & 0\\0 & \id_{E'}\end{psmallmatrix}$} (m-1-2)
			edge node [rotate=-30,below] {$\begin{psmallmatrix}g & g'\end{psmallmatrix}$} (m-2-2)
    (m-1-2) edge node [above] {$\begin{psmallmatrix}u & 0\\v & \id_{E'}\end{psmallmatrix}$} (m-1-3)
			edge node [right] {$\begin{psmallmatrix}s & g'\end{psmallmatrix}$} (m-2-2)
    (m-1-3) edge node [rotate=30,below] {$\begin{psmallmatrix}g & g'\end{psmallmatrix}$} (m-2-2);
\end{tikzpicture}
\end{center}
Since $\begin{psmallmatrix}g & g'\end{psmallmatrix}$ is right minimal, we have that
\[\begin{psmallmatrix}u & 0\\v & \id_{E'}\end{psmallmatrix} \circ \begin{psmallmatrix}t & 0\\0 & \id_{E'}\end{psmallmatrix} = \begin{psmallmatrix}u \circ t & 0\\v \circ t & \id_{E'}\end{psmallmatrix} \]
is an isomorphism.
Hence $u \circ t \colon E \to X \to E$ is an isomorphism, thus we have shown the existence of some $u \colon X \to E$ such that $u \circ t = \id_E$, so $t$ is a split monomorphism.
This finishes the proof.
\end{proof}
\end{theorem}

\begin{lemma} \label{lemma:PPDefect-2HF}
Let $A$ be a tame hereditary algebra of rank two.
Let $M$ be an indecomposable preprojective module of defect $-2$ which is not projective. Then there exists a monomorphism $f \colon N \to M$ such that $N$ is a direct sum of preprojective modules of defect $-1$ and we have that $\dim_k \coker f \leq h_A$.
\begin{proof}
First, as $M$ is an indecomposable nonprojective module, it is well-known (see e.g. \cite[Theorem~V.1.15]{AuslanderReitenSmaloe1995RepresentationTheoryArtinAlgebras}) that there exists a right minimal almost split morphism $h \colon E \to M$. As this map is the right hand side of an almost split sequence, it is surjective.
The only rank two case where this can occur is $\widetilde{A}_{11}$ of type $(1,4)$.
\proofstep{Case 1: $(1,4)$} We have $\underline{h}_A = (2,1)$.
We may assume $\underline{m} = \dimv M = (4k,2k-1)$ for some $k > 1$. Then $E = {P_{2k-1}}^{\oplus 4}$, where $P_{2k-1}$ is the unique indecomposable preprojective module of dimension vector $(2k-1,k-1)$.
Now, we pick $N$ to be the direct sum of two indecomposable summands of $E$, so $N = P_{2k-1} \oplus P_{2k-1}$. Then the induced morphism $f = h_{|N}$ is irreducible by Theorem~\ref{thm:IrreducibleFromSinkMap}. Thus $f$ is either surjective or injective. But $f$ cannot be surjective for dimension reasons, so $f$ must be injective. Moreover,
\[\dimv \coker f = \dimv M - \dimv N = (4k,2k-1)-2(2k-1,k-1) = (2,1) = \underline{h}_A.\]
Also, $\del(\coker f) = 0$.
\proofstep{Case 2: $(4,1)$} We have $\underline{h}_A = (1,2)$.
We may assume $\underline{m} = \dimv M = (2k+1,4k)$ for some $k \geq 1$. Then $E = {P_{2{k}}}^{\oplus 4}$, where $P_{2{k}}$ is the unique indecomposable preprojective module of dimension vector $(k,2k-1)$.
Now, we pick $N$ to be the direct sum of two indecomposable summands of $E$, so $N = P_{2{k}} \oplus P_{2{k}}$. Then the induced morphism $f = h_{|N}$ is irreducible by Theorem~\ref{thm:IrreducibleFromSinkMap}. Thus $f$ is either surjective or injective. But $f$ cannot be surjective for dimension reasons, so $f$ must be injective. Moreover,
\[\dimv \coker f = \dimv M - \dimv N = (2k+1,4k)-2(k,2k-1) = (1,2) = \underline{h}_A.\]
Also, $\del(\coker f) = 0$.

\afterlastproofstep Thus we have found a submodule $N$ of $M$ which is a direct sum of preprojective modules of defect $-1$ and has codimension bounded by $h_A$.
\end{proof}
\end{lemma}

\begin{proposition} \label{prop:TH2PPsHF}
Let $A$ be a finite dimensional tame hereditary $k$-algebra of rank two. Then the family of (isoclasses of) preprojective modules $\cP(A)$ is hyperfinite.
\begin{proof}
By additivity of hyperfiniteness, it is enough to check this on the indecomposable preprojective modules, which have defect either $-1$ or $-2$.
By Lemma~\ref{lemma:PPDefect-1HF}, the family of indecomposable preprojective modules of defect $-1$ is hyperfinite.
For a preprojective of defect $-2$, Lemma~\ref{lemma:PPDefect-2HF} yields the existence of a submodule, which lies in a hyperfinite family and is of bounded codimension. By Proposition~\ref{prop:ExtendingHFfromSubmodulesOfBoundedCodimension}, the indecomposable preprojectives of defect $-2$ thus form a hyperfinite family.
This concludes the proof.
\end{proof}
\end{proposition}

\begin{proposition} \label{prop:RegsHFviaPPsHF}
Let $A$ be a finite dimensional tame hereditary $k$-algebra. Assume that $\cP(A)$ is hyperfinite.
Then the family of (isoclasses of) regular modules is also hyperfinite.
\begin{proof}
By the additivity of hyperfiniteness, it is enough to show that the set of all indecomposable regular modules is a hyperfinite family.
So let $T$ be a regular indecomposable module. By Lemma~\ref{lemma:PPSubmodulesOfRegulars}, there exists a submodule $P \subset T$ of codimension bounded by the maximum of the dimensions of the indecomposable injective modules with $P$ preprojective.
Thus, an application of Proposition~\ref{prop:ExtendingHFfromSubmodulesOfBoundedCodimension} along with the hypothesis shows that $\cR(A)$ is hyperfinite.
\end{proof}
\end{proposition}

\begin{theorem} \label{thm:TH2Amenable}
Let $k$ be any field. Let $A$ be a finite dimensional tame hereditary $k$-algebra of rank two.
Then $A$ is of amenable representation type.
\begin{proof}
We will show that the preprojective, the regular and the postinjective component are each hyperfinite families to conclude the amenability of $A$. We will give an argument for the indecomposable objects in each component and then apply the additivity.

That the family of (isoclasses of) preprojective modules $\cP(A)$ is hyperfinite has been shown in Proposition~\ref{prop:TH2PPsHF}.
Moreover, the hyperfiniteness of the regular component then follows from Proposition~\ref{prop:RegsHFviaPPsHF}.

We are left to deal with the postinjective modules. By Lemma~\ref{lemma:DescentOnDefectOfPIs}, for each indecomposable postinjective module $X$, we can find a submodule $Y := \ker \theta$ (using the notation of the lemma) of strictly smaller defect. What is more, if $Y$ had a postinjective summand $Z$, it must have defect $\del(Z) < \del(X)$. We do an induction on the defect $d$.
We recursively define \[\cN_d := \cN_{d-1} \cup \{\text{indecomposable postinjectives of defect $d$}\}, \quad d \geq 1,\]
and set $\cN_0$ to be the family of all preprojective and regular modules, which is hyperfinite by additivity and the above results.
For $d=1$, the submodule $Y \subseteq X$ with $\del(X) = d$ must be in $\add \cN_0$, since there are no postinjective modules $Z$ with defect $\del(Z) < 1$. Moreover, the codimension of $Y$ is bounded by the dimension of the indecomposable injective modules, of which there are only finitely many. Hence, we can use Proposition~\ref{prop:ExtendingHFfromSubmodulesOfBoundedCodimension} to prove the hyperfiniteness of the indecomposable postinjectives of defect one. Note that the base case implies that $\cN_1$ is hyperfinite.
For the induction, note that Lemma~\ref{lemma:DescentOnDefectOfPIs} also yields a submodule in $\add \cN_d$ for every indecomposable postinjective of defect $d+1$ of bounded codimension. Assuming the hyperfiniteness of $\cN_d$, Proposition~\ref{prop:ExtendingHFfromSubmodulesOfBoundedCodimension} yields that $\cN_{d+1}$ is hyperfinite.
This finishes the proof.
\end{proof}
\end{theorem}

\subsection{Passing on amenability from subalgebras}

\begin{proposition} \label{prop:THPerpCatHF}
Let $A$ be a connected finite dimensional hereditary $k$-algebra of rank $n > 2$. 
Let all connected finite dimensional hereditary $k$-algebras of rank $n-1$ be of amenable representation type.
If $T$ is an inhomogeneous regular simple $A$-module, then $T^\perp$ is hyperfinite.
\begin{proof}
By \cite[Proposition~1.1]{GeigleLenzing1991PerpendicularCategoriesApplications}, $T^{\perp}$ is an exact abelian subcategory of $\mods A$ closed under the formation of kernels, cokernels and extensions. What is more, \cite[Theorem~4.16]{GeigleLenzing1991PerpendicularCategoriesApplications} yields that $T^{\perp} = \mods \Lambda$ for some finite dimensional hereditary algebra $\Lambda$ of rank $n-1$, along with a homological epimorphism $\varphi \colon A \to \Lambda$, which induces a functor $\varphi_* \colon \mods \Lambda \to \mods A$. 
Moreover, \cite[Theorem~10.1]{GeigleLenzing1991PerpendicularCategoriesApplications} shows that $\Lambda$ is tame and connected.

Now, if $F \colon \mods \Lambda \to T^\perp$ is an equivalence, the simple $\Lambda$-modules $S(i)$ get mapped to certain modules $B_i$ in $\mods A$.
The $k$-dimension of any module $M$ over a finite dimensional $k$-algebra is determined by the length of any composition series. Such a series for $M$ in $\Lambda$ gets mapped to a composition series in the perpendicular category, and thus a series in $\mods A$, such that the factor modules are isomorphic to some $B_i$. Letting $K' := \max_{i=1,\dots,n-1}\{\dim B_i\}$, we thus know that \[\dim_{k} {}_A F(M) \leq K' \dim_k {}_{\Lambda} M.\]
On the other hand, if $F(M) \in T^\perp$, any submodule of $F(M)$ in $T^\perp$ is also a submodule in $\mods A$, so a composition series of $F(M)$ in $\mods A$ is at least as long as one in $T^\perp$. Thus, \[\dim_k {}_{\Lambda} M \leq \dim_k {}_A F(M),\] using the fact that the length of $M$ in $\mods \Lambda$ equals the length of $F(M)$ considered as an object of $T^\perp$.
Hence by Proposition~\ref{prop:HFPreservingFunctors}, we have that each $T^{\perp}$ is a hyperfinite family.
\end{proof}
\end{proposition}

We will introduce some notation for the indecomposable modules of a given homogeneous tube $\T$ of a finite dimensional tame hereditary algebra $A$. We  start by denoting the isoclasses of regular simples on the mouth of $\T$ by $T_1, \dots, T_m$ such that $\tau T_i = T_{i-1}$ for $i = 2, \dots, m$ and $\tau T_1 = T_m$. Here, $m$ is the rank of $\T$.
Similar to \cite[Chapter~3]{Ringel1984TameAlgebrasIntegralQuadraticForms}, we then define the objects $T_i[\ell]$. First, let $T_i[1] := T_i$ for each $1 \leq i \leq m$. Now, for $\ell \geq 2$, recursively define 
$T_i[\ell]$ to be the indecomposable module in $\T$ with $T_i[1]$ as a submodule such that $T_i[\ell] / T_i[1] \isom T_{i+1}[\ell-1]$. Thus $T_i[\ell]$ is the regular module of regular length $\ell$ with regular socle $T_i$. Any regular indecomposable in $\T$ will be given as some $T_i[\ell]$ since $\T$ is uniserial. We may define $T_i[\ell]$ for all $i \in \Z$ by letting $T_i[\ell] \isom T_j[\ell]$ iff $i \equiv j \mod m$.
Note that \begin{equation} \label{eq:DimvIndecReg} \dimv T_i[\ell] = \sum_{j=i}^{i+\ell-1} \dimv T_j.\end{equation}

\begin{lemma} \label{lemma:SumRegularSimplesIshA}
Let $\T$ be an inhomogeneous tube of rank $m$. 
Then $\sum_{i=1}^{m} \dimv T_i = g_\T h_A$, and $g_\T$ is globally bounded by some $g$. 
\begin{proof}
Given an inhomogeneous tube of rank $m \geq 2$ with regular simples $T_1,\dots,T_m$, recall that e.g. $T_1[m]$ is homogeneous, i.e. $\dimv T_1[m]$ is a multiple of $\underline{h}_A$.
By \eqref{eq:DimvIndecReg}, this implies that the sum of the dimension vectors of the regular simples in each tube is a multiple of $\underline{h}_A$. As there are only finitely many inhomogeneous tubes, there is a global bound on this multiple. 
\end{proof}
\end{lemma}

\begin{lemma} \label{lemma:THRegularsPerps}
Let $\T$ be a tube of rank $m \geq 2$. 
Let $X$ be an indecomposable regular in $\T$.
Then there exists a submodule $Y \subseteq X$ of codimension bounded by the sum of the $k$-valued sum of the entries of $g \underline{h}_A$ and a regular simple module $T \in \T$ such that $Y \in T^{\perp}$.
\begin{proof}
By adapting the proof of \cite[3.1.3']{Ringel1984TameAlgebrasIntegralQuadraticForms} to the fact that $\dim_k \Hom(T_i,T_i) = e$ for all $i$, we have that
\[ \langle \dimv T_i, \dimv T_j \rangle =
\begin{cases}
\phantom{-}e, & i \equiv j\phantom{+1} \mod m, \\
-e, & i \equiv j+1 \mod m, \\
 \phantom{-}0, & \text{else.}
\end{cases} \] 
This implies that for any $j \in \Z$, $\langle \dimv T_j, \dimv T_i[\ell]\rangle = 0$, provided $\ell \equiv 0 \mod m$. For $T_i[\ell]$ is uniserial, we also have that $\Hom(T_j,T_i[\ell]) = 0$ if and only if $j \not\equiv i \mod m$. 

We write $\ell = n \cdot m + r$, where $0 \leq r < m$. Then there is a short exact sequence
\[ 0 \to T_i[nm] \to T_i[\ell] \to Z \to 0,\]
where $\dimv Z \leq g \underline{h_A}$ using \eqref{eq:DimvIndecReg} and Lemma~\ref{lemma:SumRegularSimplesIshA}. Thus, we have found a suitable submodule $T_i[nm] \in T_{i+1}^\perp$.
\end{proof}
\end{lemma}

\begin{lemma} \label{lemma:PerpForPP}
Let $A$ be a finite dimensional tame hereditary algebra. Let $S$ be some regular simple inhomogeneous module in a tube $\T$, then there is a constant $m >0$. 
If $X=\tau^{-p} P(i)$ is some indecomposbale preprojective $A$-module, there exists a submodule $Y \subseteq X$ and a module $Q$ with $\dimv Q \leq m \dimv S$ such that $0 \to Y \to X \to Q \to 0$ is exact and $Y \in {\tau^{-}S}^{\perp}$.
\begin{proof}
Let $f \colon X \to S^m$ be a minimal left $\add(S)$-approximation of $X$, in particular, the induced map \[f^* \colon \Hom(S^m,S) \to \Hom(X,S), \quad \phi \mapsto \phi \circ f\]
is surjective. 
Now, let $Y = \ker f$ and $Q = \im f$. Note that $f = g \circ h$ factors through $h \colon X \to Q$ and $g \colon Q \to S^m$. We have a short exact sequence
\[\xi \colon 0 \to Y \xrightarrow{\iota} X \xrightarrow{h} Q \to 0,\]
to which we apply $\Hom(-,S)$ and get the long exact sequence
\begin{multline*} 0 \to  \Hom(Q,S) \xrightarrow{h^*} \Hom(X,S) \xrightarrow{\iota^*} \Hom(Y,S) \xrightarrow{\xi^*} \Ext{}{1}(Q,S)\\ \to \Ext{}{1}(X,S) \to \Ext{}{1}(Y,S) \to 0,\end{multline*}
for $A$ is hereditary.
First, by \cite[2.4.6*]{Ringel1984TameAlgebrasIntegralQuadraticForms}, $\Ext{}{1}(X,S) \isom D\Hom(\tau S,X) = 0$ since there are no maps from the regular to the preprojective component.
Second, $h^*$ is injective. But $f^* = (g \circ h)^* = h^* \circ g^*$ is surjective, so $h^*$ is also surjective, hence bijective.
Since $Q$ can only have preprojective summands or regular summands isomorphic to copies of $S$, we must have $\Ext{}{1}(Q,S) = 0$, implying that $\Hom(Y,S) = 0$.
This shows that $Y \in {^\perp} S = (\tau^{-}S)^\perp$.

By the construction of the (minimal) left $\add(S)$-approximation of $X$, we see that \[m = \frac{\dim_k \Hom(X,S)}{\dim_k \End(S)}.\]
Denoting the regular simples on the mouth of $\T$ by $T_1,\dots,T_m$, we may assume that $S \isom T_j$ for some $1 \leq j \leq n$.
Then, \begin{align*} \dim_k \Hom(X,S) - \underbrace{\dim_k \Ext{}{1}(X,S)}_{=0} &= \langle \dimv X, \dimv T_{j} \rangle = \langle \dimv \tau^{-p} P(i), \dimv T_{{j}} \rangle \\
& = \langle \dimv P(i), \dimv T_{{j-p}} \rangle = f_i (\dimv T_{j-p})_i,\end{align*}
where $f_i = \dim_k \End(P(i))$.
Now, $\dimv Q \leq \dimv S^m \leq m \dimv S$.
Thus, the codimension of $Y$ as a submodule of $X$ is bounded by a constant solely depending on the finitely many regular simple inhomogeneous modules and their finite $k$-dimensions.

\end{proof}
\end{lemma}

\begin{theorem} \label{thm:THAmenable}
Let $k$ be any field. Let $A$ be a finite dimensional tame hereditary $k$-algebra. Then $A$ is of amenable representation type.
\begin{proof}
We do an induction on the rank of $A$. There are no non-finite tame hereditary $k$-algebras of rank one.
The case $n=2$ is the subject of Theorem~\ref{thm:TH2Amenable}.

Let $A$ be some finite dimensional tame hereditary algebra of rank $n>2$. Now assume that it has been shown that all finite dimensional tame hereditary connected $k$-algebras of rank $n-1$ are of amenable type.

We will show that the indecomposable preprojective modules form a hyperfinite family first.
By an application of Lemma~\ref{lemma:PerpForPP}, all indecomposable preprojective modules $P$ have a submodule of (globally) bounded codimension which lies in the perpendicular category of a fixed regular simple module  $T_1$ in an inhomogeneous tube $\T$. Then Proposition~\ref{prop:THPerpCatHF} shows that the indecomposable preprojectives lie in a hyperfinite family, additionally applying Proposition~\ref{prop:ExtendingHFfromSubmodulesOfBoundedCodimension}. This uses the fact that there are only finitely many inhomogeneous tubes, each of finite rank.

Next, we consider the regular modules. Indecomposable regular modules in a tube other than $\T$ will be contained in $T_1^{\perp}$ by an analogue of \cite[3.1.3']{Ringel1984TameAlgebrasIntegralQuadraticForms}.
By Lemma~\ref{lemma:THRegularsPerps}, any regular indecomposable in $\T$ either is contained in the perpendicular category of some regular simple in $\T$ or has a submodule of bounded codimension, which is in the perpendicular category of some regular simple in $\T$. But by the Proposition~\ref{prop:THPerpCatHF}, the perpendicular categories are hyperfinite. In the latter case, we can therefore apply Proposition~\ref{prop:ExtendingHFfromSubmodulesOfBoundedCodimension} to show the hyperfiniteness of the family of these indecomposable regular modules.

For the postinjective modules, we may apply the argument used for the postinjectives in the proof of Thereom~\ref{thm:TH2Amenable}. This completes the induction and finishes the proof.
\end{proof}
\end{theorem}

\section{Tilted Algebras}

Let $A = kQ$ be the finite dimensional path algebra of a finite acyclic quiver. 

Recall that a module $T$ is a \emph{tilting module} for $A$ provided 
\begin{itemize}
\item $T$ is rigid, i.e. $\Ext{}{1}(T,T) = 0$,
\item $\pdim T \leq 1$, (always true for $A$-modules),
\item the number of isoclasses of indecomposable direct summands of $T$ is equal to the number of vertices of $Q$ (resp. the number of isoclasses of simple modules)
\end{itemize}

Now, an algebra of the form $B = \End_A(T)$ for some tilting module $T$ is called a \emph{tilted algebra}.

\begin{theorem}
Let $A$ be a finite dimensional tame hereditary algebra. Let $T$ be a tilting module without a postinjective summand. Then $\End_A(T)$ is of amenable representation type.
\begin{proof}
First note that $T$ cannot be regular by \cite[Lemma~3.1]{HappelRingel1981ConstructionTiltedAlgebras}, since $A$ is tame. Thus $T$ must have a nonzero preprojective summand.
Now, by \cite[Proposition~3.2]{HappelRingel1981ConstructionTiltedAlgebras}, the torsion class $\cT(T)$ is infinite. It is in bijection with the torsion free class $\cY(T)$ via $G = \Hom(T,-)$ of the Brenner-Butler tilting theorem. But $\cY(T)$ is hyperfinite, since $G$ is an HF-preserving functor using Proposition~\ref{prop:HFPreservingFunctors} and the amenability of $\mods A$ from Theorem~\ref{thm:THAmenable}.
Since $T$  is a splitting tilting module (see e.g. \cite[Corollary~VI.5.7]{AssemSimsonSkowronski2010TechniquesRT}), any indecomposable module of $\mods B$ either lies in $\cY(T)$ or $\cX(T)$, the torsion class in $\mods B$. But $\cF(T)$, which is in bijection to $\cX(T)$ via the Tilting Theorem, is finite by \cite[Proposition~3.2*]{HappelRingel1981ConstructionTiltedAlgebras}\footnote{note that condition (vii) in the cited proposition should read ``$T_A$ has no non-zero preprojective direct summand.''}.
\end{proof}
\end{theorem}

\subsection{Concealed algebras}
In this subsection, let $k = \overline{k}$ denote an algebraically closed field.

Recall that a module $M \in \mods A$ is preprojective if each indecomposable summand of $M$ lies in some preprojective component of $\Gamma_A$, the AR-quiver of $A$.

Now, let $A$ be hereditary and $T \in \mods A$ be some preprojective tilting module. Then the tilted algebra $B = \End_A(T)$ is called a \emph{concealed algebra}.

\begin{corollary} \label{cor:TCAmenable}
Let $B = \End_A(T)$ be a tame concealed algebra, i.e. let the algebra $A$ be of tame but not finite representation type.
Then $B$ is of amenable representation type.
\end{corollary}

\begin{acknowledgement}
These notes are based on work done during the author’s doctorate studies at Bielefeld University. The author would like to thank his supervisor Professor W. Crawley-Boevey for his advice and guidance and A.~Hubery for helpful discussions and a suggestion for a nicer proof of Lemma~\ref{lemma:PerpForPP}.
\end{acknowledgement}

{
 \printbibliography
}
\end{document}